\documentclass[12pt]{article}
\usepackage{amsfonts, amssymb, dsfont}
\usepackage{graphicx}
\usepackage{epstopdf}
\usepackage{algorithmic}
\usepackage{color}
\ifpdf
  \DeclareGraphicsExtensions{.eps,.pdf,.png,.jpg}
\else
  \DeclareGraphicsExtensions{.eps}
\fi

\usepackage{amsopn}
\usepackage{mathtools}
\usepackage{subfigure}
\usepackage[left=0.9in,right=0.9in,top=1in,bottom=1in,includefoot,includehead,headheight=13.6pt]{geometry}

\newcommand{\cD}{\ensuremath{\mathcal{D}}}

\newcommand{\cM}{\ensuremath{\mathcal{M}}}

\newcommand{\cP}{\ensuremath{\mathcal{P}}}

\newcommand{\cR}{\ensuremath{\mathcal{R}}}

\newcommand{\bC}{\ensuremath{\mathbb{C}}}

\newcommand{\bE}{\ensuremath{\mathbb{E}}}

\newcommand{\bP}{\ensuremath{\mathbb{P}}}

\newcommand{\bR}{\ensuremath{\mathbb{R}}}
\newcommand{\bS}{\ensuremath{\mathbb{S}}}

\def\[{\left[}
\def\]{\right]}
\def\<{\langle}
\def\>{\rangle}
\def\({\left(}
\def\){\right)}
\def\[{\left [}
\def\]{\right]}
\def\({\left(}
\def\){\right)}

\newcommand{\Vn}{\ensuremath{{V_n}}}
\newcommand{\Wm}{\ensuremath{{W_m}}}

\newcommand{\wc}{\ensuremath{\mathrm{wc}}}
\newcommand{\ms}{\ensuremath{\mathrm{ms}}}

\newcommand{\pbdw}{\ensuremath{\text{(pbdw)}}}
\newcommand{\wcpbdw}{\ensuremath{\text{(wc, pbdw)}}}
\newcommand{\mspbdw}{\ensuremath{\text{(ms, pbdw)}}}

\newcommand{\vphi}{\ensuremath{\phi}}
\newcommand{\keff}{\ensuremath{k_{\text{eff}}}}

\newcommand{\dist}{\operatorname{dist}}

\DeclareMathOperator*{\argmin}{arg\,min}

\newcommand{\rd}{\ensuremath{\mathrm d}}

\newcommand{\eps}{\ensuremath{\varepsilon}}

\newcommand{\cond}{\; :\;}

\usepackage{soul}

\newcommand{\transport}{\ensuremath{\text{tr}}}
\newcommand{\diffusion}{\ensuremath{\text{diff}}}
\usepackage{tikz} 
\usepackage[colorlinks = true,
linkcolor = blue,
urlcolor  = blue,
citecolor = blue,
anchorcolor = blue]{hyperref}

\title{Impact of physical model error on State Estimation for neutronics applications}
\date{}
\author{Y. Conjungo Taumhas, D. Labeurthre, F. Madiot, O.~Mula, T.~Taddei}

\begin{document}

\maketitle

\begin{abstract}
In this paper, we consider the inverse problem of  state estimation of nuclear power fields in a power plant from a limited number of observations of the neutron flux. For this, we use the Parametrized Background Data Weak approach. The method combines the observations with a parametrized PDE model for the behavior of the neutron flux. Since, in general, even the most sophisticated models cannot perfectly capture reality, an inevitable model error is made. We investigate the impact of the model error in the power reconstruction when we use a diffusion model for the neutron flux, and assume that the true physics are governed by a neutron transport model.
\end{abstract}
%
%
%
\section*{Introduction }

In the field of nuclear engineering, numerical methods play a crucial role at several stages: they are involved in important assessments and decisions related to design, safety, energy efficiency, and reactor loading plans. In this paper, we focus on the task of providing real time information about the spatial distribution of the nuclear power generated by a nuclear reactor from a limited number of measurement observations. We combine this data with physical models in order to provide a complete spatial reconstruction of the power field. This task is a state estimation problem, and we work with the Parametrized Background Data Weak (PBDW), originally introduced in \cite{MPPY2015}. The method has the appealing feature of providing very fast reconstructions by leveraging techniques from model order reduction of parametric Partial Differential Equations (PDEs). We refer to \cite{binev2017, DPW2017, CDDFMN2020, CDMN2022} for theoretical analysis of the method, optimal recovery results and nonlinear extensions. A recent overview may be found in~\cite{Mula2022}.


The main ideas of the above state estimation methodology have been applied to the field of nuclear physics for applications connected to neutronics (see \cite{ABGMM2017, ABGMM2017-proc, ABCGMM2018}). We could also cite other works such as \cite{SCLR2016,GCZC2020,ChenEtAll2021,GERMAN2019144,LORENZI2018245} which study the forward reduced modeling problem for neutronics (compared to these works, note that there is a salient difference in the nature of the task that we consider, which is inverse state estimation). In this paper, we again consider neutronics but our goal is to study the impact of inaccuracies in the physical model that is involved in the reconstruction algorithm, and which is often assumed to perfectly describe reality. This assumption goes beyond the present application on neutronics but studying it for this particular topic has the advantage that we have two very well identified models with different levels of accuracy, thereby allowing to examine synthetically what one can expect when working in a real application scenario.

In neutronics, the most accurate physical model is the so-called neutron  transport equation which describes the evolution of the neutronic population in a reactor core by expressing it in the form of a balance between produced and lost neutrons~\cite{Hebert2009}. This model is often approximated at the reactor core scale by a neutron diffusion model to save computing time. This is why in this work, we explore the impact of model inaccuracies by applying a reconstruction based on a diffusion model for the neutron flux, and then  assuming that the true physical system is  governed by a neutron transport model.

The paper is organized as follows. Section~\ref{sec:inverse_state_estimation} is devoted to presenting inverse state estimation problems and the PBDW method. Section~\ref{sec:application_reconstruction_power} details the application of the methodology to the reconstruction of nuclear power. Section~\ref{sec:numerics} provides some numerical results.

%
%
%
%
%

\section{Inverse State Estimation with PBDW}
\label{sec:inverse_state_estimation}
In this section, we introduce the problem of state estimation, and the Parametrized Background Data Weak method which combines measurement observations and reduced models from parametric PDEs. We refer the reader to \cite{Mula2022} for an overview of inverse problem algorithms using these elements.

Let $\cR$ be a fixed given domain of $\bR^d$ with dimension $d\geq 1$, and let $V$ be a Hilbert space defined over $\cR$. In our application, $\cR$ will be defined as the nuclear reactor domain. The space is endowed with an inner product $\<\cdot, \cdot\>$ and induced norm $\Vert \cdot \Vert$. The choice of $V$ must be relevant for the problem under consideration:
 typical options are $L^2$, $H^1$; for pointwise measurements, a Reproducing Kernel Hilbert Space should be considered.

Our goal is to recover an unknown function $u\in V$ from $m$ measurement observations
\begin{equation}
\label{eq:meas}
y_i = \ell_i(u),\quad i=1,\dots,m,
\end{equation}
where the $\ell_i$ are linearly independent linear forms from the dual $V'$. Note that  
we have assumed that experimental observations are perfect; however, the methodology could be extended to deal with noisy measurements (see, e.g., \cite{Taddei2017, GMMT2019, EF2021}). In practical applications, each $\ell_i$ models a sensor device which is used to collect the measurement data $\ell_i(u)$. In the applications that we present in our numerical tests, the observations come from sensors for the neutron flux which are placed in the reactor.

We denote by $\omega_i \in V$ the Riesz representers of the $\ell_i$. They are defined via the variational equation
$$
\left< \omega_i, v \right> = \ell_i(v),\quad \forall v \in V.
$$
Since the $\ell_i$ are linearly independent in $V'$, so are the $\omega_i$ in $V$ and they span an $m$-dimensional space
$$
\Wm={\rm span}\{\omega_1,\dots,\omega_m\} \subset V.
$$

When there is no measurement noise, knowing the observations $y_i=\ell_i(u)$ is equivalent to knowing the orthogonal projection
\begin{equation}
\omega =P_{\Wm} u.
\end{equation}
In this setting, the task of recovering $u$ from the measurement observation $\omega$ can be viewed as building a recovery algorithm
$$
A:\Wm\mapsto V
$$
such that $A(P_{\Wm}u)$ is a good approximation of $u$ in the sense that $\Vert u - A(P_{\Wm}u) \Vert$ is small.

Recovering $u$ from the measurements $P_{\Wm} u$ is a very ill-posed problem since $V$ is generally a space of very high or infinite dimension so, in general, there are infinitely many $v\in V$ such that $P_{\Wm} v = \omega$. It is thus necessary to add some a priori information on $u$ in order to recover the state up to a guaranteed accuracy. In the following, we work in the setting where $u$ is a solution to some parameter-dependent PDE of the general form 
\begin{equation*}
\cP(u, \mu) = 0,
\end{equation*}
where $\cP$ is a differential operator and $\mu$ is a vector of parameters that describe some physical property and  belong to a  given set $\cD\subset \bR^p$. For every $\mu \in \cD$, we assume that the PDE has a unique solution $u=u(\mu)\in V$. Therefore, our prior on $u$ is that it belongs to the so-called {\em solution manifold}
\begin{equation}
\label{eq:manifold}
\cM \coloneqq \{ u(\mu) \in V \; :\; \mu\in\cD \}.
\end{equation}
In practical applications, the PDE model $\mathcal{P}$ might not be known exactly or might be too expensive to evaluate:  we should thus rely on a surrogate approximate model to perform state estimation.

\medskip
\noindent
\textbf{Performance Benchmarks:} The quality of a recovery mapping $A$ is quantified in two ways:
\begin{itemize}
\item If the sole prior information is that $u$ belongs to the manifold $\cM$, the performance is usually measured by the worst case reconstruction error
\begin{equation}
\label{eq:err-A-wc}
E_{\wc}(A,\cM) = \sup_{u\in\cM} \Vert u - A(P_{\Wm} u) \Vert \, .
\end{equation}
\item In some cases $u$ is described by a probability distribution $p$ on $V$ supported on $\cM$. This distribution is itself induced by a probability distribution on $\cD$ that is assumed to be known. When no information about the distribution is available, usually the uniform distribution is taken. In this Bayesian-type setting, the performance is usually measured in an average sense through the mean-square error
\begin{equation}
\label{eq:err-A-ms}
E^2_{\ms}(A,\cM) = \bE\left( \Vert u - A(P_{\Wm} u) \Vert^2\right) = \int_{V} \Vert u - A(P_{\Wm} u)\Vert^2 dp(u) \, ,
\end{equation}
and it naturally follows that $E_{\ms}(A,\cM)\leq E_{\wc}(A,\cM) $.
\end{itemize}

\medskip
\noindent
\textbf{PBDW algorithm:} In this work, we   resort to the Parametrized-Background Data-Weak algorithm (PBDW, \cite{MPPY2015}) to estimate the state $u$. Other choices would of course be possible but the PBDW algorithm is relevant for the following reasons:
\begin{itemize}
\item \textbf{Simplicity and Speed:} It is easily implementable and it provides reconstructions in near-real time.
\item \textbf{Optimality:} It has strong connections with optimal linear reconstruction algorithms as has been studied in \cite{BCDDPW2017, CDDFMN2020}.
\item \textbf{Extensions:} If required, the algorithm can easily be extended to enhance its reconstruction performance (see \cite{CDMN2022, GGLM2021}). In particular, it is shown in \cite{CDMN2022} that piecewise PBDW reconstruction strategy can deliver near-optimal performance. The PBDW algorithm can also be easily adapted to accommodate noisy measurements (see \cite{Taddei2017, GMMT2019}) and some easy-to implement extension to mitigate the model error exist (in the following however, we assume the PDE model is perfect for the sake of simplicity).
\end{itemize}
Since the geometry of $\cM$ is generally complex, optimization tasks posed on $\cM$ are difficult (lack of convexity, high evaluation costs for different parameters). Therefore, instead of working with $\cM$, PBDW works with a linear (or affine) space $V_n$ of reduced dimension $n$ which is expected to approximate the solution manifold well in the sense that the approximation error of the manifold
\begin{equation}
\label{eq:error-manifold}
\delta_n^{(\wc)} \coloneqq \sup_{u\in\cM} \dist(u, V_n)  \,
,\quad \text{or} \quad
\delta^{(\ms)}_n \coloneqq \bE\left(  \dist(u, V_n)^2\right)^{1/2} \, 
\end{equation}
decays rapidly if we increase the dimension $n$. It has been proven in \cite{CD2015acta} that it is possible to find such hierarchies of spaces $(V_n)_{n\geq 1}$ for certain manifolds coming from classes of elliptic and parabolic problems, and numerous strategies have been proposed to build the spaces in practice (see, e.g., \cite{BMPPT2012,RHP2007} for reduced basis techniques and \cite{CD2015acta,CDS2011} for polynomial approximations in the $\mu$ variable).

Assuming that we are given a reduced model $V_n$ with $1\leq n\leq m$, the PBDW algorithm $$A^\pbdw_{m,n}:\Wm \to V$$ gives for any $\omega \in \Wm$ a solution of
\begin{equation}
\label{eq:pbdw-algo}
A^\pbdw_{m,n}(\omega) \in \argmin_{u \in \omega + W(\cR)^\perp } \dist(u, V_n).
\end{equation}
The minimizer is unique as soon as $n\leq m$ and $\beta(V_n,\Wm)>0$, which is an assumption to which we adhere in the following. The quantity $\beta$ is defined as follows. For any pair of closed subspaces $(E,F)$ of $V$, $\beta(E, F)$ is defined as
\begin{equation}
\label{eq:infsup}
\beta(E,F):=\inf_{e\in E}\sup_{f\in F}\frac {\<e,f\>}{\|e\|\, \|f\|}=\inf_{e\in E}\frac {\|P_{F} e\|}{\|e\|} \in [0,1] .
\end{equation}
We can prove that $A^\pbdw_{m,n}$ is a bounded linear map from $\Wm$ to $V_n\oplus(W_m\cap V_n^\perp)$. 

In practice, solving problem \eqref{eq:pbdw-algo} boils down to solving a linear least squares minimization problem whose cost is essentially of order $n^2+m$, and we can compute $\beta(V_n, W_m)$ by finding the smallest eigenvalue of an $n\times n$ matrix. We refer, e.g., to \cite[Appendix A, B]{Mula2022} for details on how to compute these elements in practice. It follows that, since in general $m$ is not very large, if the dimension $n$ of the reduced model is moderate, the reconstruction with \eqref{eq:pbdw-algo} can take place in close to real-time.

For any $u\in V$, the reconstruction error is bounded by
\begin{equation}
\Vert u - A^\pbdw_{m,n}(\omega) \Vert
\leq \beta^{-1}(V_n, \Wm) \Vert u - P_{\Vn\oplus(W_m\cap V_n^\perp) } u \Vert
\leq \beta^{-1}(V_n, \Wm) \Vert u - P_{\Vn}u \Vert,
\label{eq:pbdw_bound}
\end{equation}
where we have omitted the dependency of the spaces on $\cR$ in order not to overload the notation, and we will keep omitting this dependency until the end of this section. Depending on whether $\Vn$ is built to address the worst case or mean square error, the reconstruction performance over the whole manifold $\cM$ is bounded by
\begin{equation}
\label{eq:err-wc-pbdw}
e_{m,n}^\wcpbdw \coloneqq
E_{\wc}(A^\pbdw_{m,n}, \cM) \leq \beta^{-1}(V_n, \Wm) \max_{u\in \cM}\dist(u, V_n\oplus(V_n^\perp\cap\Wm)) \leq \beta^{-1}(V_n, \Wm) \, \delta_n^{(\wc)},
\end{equation}
or
\begin{align}
e_{m,n}^\mspbdw \coloneqq
E_{\ms}(A^\pbdw_{m,n}, \cM)
&\leq \beta^{-1}(V_n, \Wm) \bE\left(\dist(u, V_n\oplus(V_n^\perp\cap\Wm))^2\right)^{1/2}  \nonumber\\
&\leq \beta^{-1}(V_n, \Wm) \, \delta_n^{(\ms)} . \label{eq:err-ms-pbdw}
\end{align}
Note that $\beta(\Vn, \Wm)$ can be understood as a stability constant. It can also be interpreted as the cosine of the angle between $\Vn$ and $\Wm$. The error bounds involve the distance of $u$ to the space $V_n\oplus(V_n^\perp\cap\Wm)$ which provides slightly more accuracy than the reduced model $\Vn$ alone. This term is the reason why it is sometimes said that the method can correct model error to some extent. In the following, to ease the reading we will write errors only with the second type of bounds \eqref{eq:err-ms-pbdw} that do not involve the correction part on $V_n^\perp\cap\Wm$.

An important observation is that for a fixed measurement space $W_m$ (which is the setting in our numerical tests), the error functions
$$
n\mapsto e_{m,n}^\wcpbdw,\quad \text{and} \quad n\mapsto e_{m,n}^\mspbdw
$$
reach a minimal value for a certain dimension $n^*_\wc$ and $n^*_\ms$ as the dimension $n$ varies from 1 to $m$. This behavior is due to the trade-off between:
\begin{itemize}
\item the improvement of the approximation properties of $V_n$ as $n$ grows ($\delta_n^{(\wc)}$ and $\delta_n^{(\ms)} \to 0$ as $n$ grows)
\item the degradation of the stability of the algorithm, given here by the decrease of $\beta(V_n, \Wm)$ to 0 as $n\to m$. When $n> m$, $\beta(V_n, \Wm)=0$.
\end{itemize}
As a result, the best reconstruction performance with  PBDW is given by
$$
e_{m,n^*_\wc}^\wcpbdw = \min_{1\leq n \leq m} e_{m,n}^\wcpbdw,
\quad\text{or}\quad
e_{m,n^*_\ms}^\mspbdw = \min_{1\leq n \leq m} e_{m,n}^\mspbdw.
$$
\paragraph{{\textbf{Noise and Model Error:}}} To account for measurement noise and model bias in the above analysis, let us assume that we get noisy observations $\tilde \omega = \omega + \eta$ with $||\eta|| \leq \eps_{noise}$. Suppose also that the true state $u$ does not lie in $\cM$ but satisfies $\dist(u, \cM) \leq \eps_{model}$. We can prove that the error bound \eqref{eq:pbdw_bound} should be modified into
\begin{align*}
\Vert u - A^\pbdw_{m,n}(\tilde \omega) \Vert
&\leq \beta^{-1}(V_n, \Wm) (\Vert u - P_{\Vn}u \Vert + \eps_{noise} + \eps_{model}).
\end{align*}
Thus \eqref{eq:err-wc-pbdw} and \eqref{eq:err-ms-pbdw} become
\begin{equation}
e_{m,n}^\wcpbdw \coloneqq
E_{\wc}(A^\pbdw_{m,n}, \cM)
\leq \beta^{-1}(V_n, \Wm) \, (\delta_n^{(\wc)}+ \eps_{noise} + \eps_{model}),
\end{equation}
and
\begin{equation}
e_{m,n}^\mspbdw \coloneqq
E_{\ms}(A^\pbdw_{m,n}, \cM)
\leq \beta^{-1}(V_n, \Wm) \, (\delta_n^{(\ms)}+ \eps_{noise} + \eps_{model}).
\end{equation}
Note that the estimation accuracy benefits from decreasing the model error, and the noise. Since both errors have the same additive effect on the reconstruction accuracy, model error could be understood as measurement error and vice-versa.  However, since the underlying physical reasons leading to model and measurement error are entirely different, it is preferable to clearly keep both concepts separately. Note further that the computational complexity of the method is not affected by these errors. This is in contrast to Bayesian methods for which small noise levels induce computational difficulties due to the concentration of the posterior distribution. 

\paragraph{{\textbf{Sensor modeling error:}}} Another error that can occur comes from our choice of the observation functions $\omega_i$ which are built to mimic the response of the sensor devices. Suppose that we work with imperfect functions $\tilde \omega_i$ that deviate from the exact one $\omega_i$ with $\Vert \omega_i - \tilde \omega_i \Vert \leq \rho$ for some $\rho>0$. Then noiseless observations can be written as
$$
y_i = \ell_i(u) = \left< \omega_i, u \right> = \left< \tilde \omega_i, u\right> + \left< \omega_i - \tilde \omega_i, u\right>.
$$
The right hand side tells us that by working with the inexact $\tilde \omega_i$, we are introducing a term of noise which is $ \left< \omega_i - \tilde \omega_i, u\right>$. The noise has level $\rho \Vert u \Vert$. It follows that working with an inexact representation of the sensor response can be understood as introducing additional noise to the observations.

\section{Application to the reconstruction of nuclear power}
\label{sec:application_reconstruction_power}

In this work, we apply the above general framework to reconstruct the nuclear power $P$ generated in a nuclear reactor core defined on a convex domain $\cR$. The power
$P$ is a real-valued function in $\cR$, $P:\cR\to \bR_+$, and in the following we  reconstruct it by viewing it as a function in the space
$$
V= L^2(\cR).
$$

The nuclear power $P$ we want to rebuild always comes from the neutron transport model. However, the spaces used to reconstruct $P$ will be divided in two cases. One space is made up of solutions of the transport model while the other is made up of solutions of the diffusion model as discussed in the following sections.

\subsection{The neutron transport model}
We assume that the reactor is in a stationary state where the neutron population  $\psi$, usually called the angular flux, depends on $(r, \omega, E)$, namely the spatial position $r \in \cR\subset \bR^d$, the direction of propagation $\omega\in \mathbb{S}_{d}$ where $\mathbb{S}_{d}$ is the unit sphere of $\bR^d$, and the kinetic energy $E\in \bR^+$. We work with a multi-group approach where we consider a discrete set of energies $E_G< \dots < E_{0}$, and we denote
$$
\psi(r, \omega, [E_{g},E_{g-1}]) \coloneqq 
\psi^g(r, \omega),
\quad \forall (r,\omega)\in \cR\times \bS_d,
\quad \forall g\in \{1,\dots, G\}.
$$

With this notation, the neutron transport equation is a generalized eigenvalue problem in which we search for a multigroup flux $\psi = (\psi^g)_{g=1}^G$, and a generalized eigenvalue $\lambda\in \bC^*$ (see~\cite{DLvol6})




\begin{equation}
\begin{cases}  
L^g\psi^g(r,\omega) &= H^g \psi(r,\omega)+ \lambda F^g\psi^g(r,\omega) \quad \text{ in } \cR\times\mathbb{S}_{2},  \quad \forall g\in \{1,\dots, G\}
\\
\psi(r, \omega) &=0  \text{\hspace{4.cm} on } \partial \Gamma_- \coloneqq\{(r, \omega) \in \partial \cR \times \bS_d\cond n(r)\cdot \omega<0\},
\end{cases}
\label{eq:boltzmann_void}
\end{equation}
where 
\begin{align*}
&L^g\psi^g(r,\omega) \coloneqq (\omega \cdot \nabla + \Sigma_t^ g(r))\psi^g(r,\omega) \text{ is the advection operator},\\
&H^g\psi(r,\omega)\coloneqq \sum_{g'=1}^G\int_{ \mathbb{S}_2}\Sigma_s^{g'\rightarrow g}(r,\omega'\cdot\omega) \psi^{g'}(r,\omega')d\omega'\text{ is the scattering operator},\\
&F^g\psi(r,\omega)\coloneqq\frac{\chi^g(r)}{4\pi}\sum_{g'=1}^G (\nu\Sigma_f)^{g'}(r) \int_{\mathbb{S}_2} \psi^{g'}(r,\omega)d\omega  \text{ is the fission operator}.
\end{align*}
In the listed terms, $\Sigma_t^g(r)$ denotes the total cross-section and $\Sigma_s^{g'\rightarrow g}(r,\omega'\cdot\omega) $ is the scattering cross-section from energy group $g'$ and direction $\omega'$ to energy group $g$ and direction $\omega$, $\Sigma_f^g(r)$ is the fission cross-section, $\nu^g( r)$ is the average number of neutrons emitted per fission and $\chi^g(r)$ is the fission spectrum. We suppose that all the coefficients are measurable bounded functions of their arguments.

Under certain conditions (which we assume to be satisfied in the following), the eigenvalue $\lambda_{\min}$ with the smallest modulus is simple, real and strictly positive. We refer to~\cite[Theorem 2.2]{allaire1999homogenization} for the sketch of the proof detailed in~\cite[Theorem 2.1.1, p 92]{balThesis}. The associated eigenfunction $\psi$ belongs to the Hilbert space $W^2(\cR)^G$ where $W^2(\cR\times \mathbb{S}_2)=\{\psi\in L^2(\cR\times \mathbb{S}_2) \text{ s.t. } \omega\cdot \nabla \psi \in L^2(\cR\times \mathbb{S}_2)\}$, is also real and positive at almost every $(x,\omega)\in\cR \times \mathbb{S}_2$. With this model, once the neutron flux is computed by solving \eqref{eq:boltzmann_void} numerically, the nuclear power is given by
$$
P({r}) \coloneqq \sum_{g'=1}^G  (\kappa \Sigma_f)^{g'} {\int_{\bS_2}}\psi^{g'}({r}, \omega) \,\rd\omega,\quad \forall {r}\in \cR \text{ a.e,}
$$
where $\kappa^g\in L^{\infty}(\cR)$ is the released energy per fission and since $\psi \in (W^2(\cR\times \mathbb{S}_2))^G$, we have that $P\in V$.

\subsection{The neutron diffusion equations}
In this work, the neutron flux $\phi$ is modeled with the two-group neutron diffusion equation with null flux boundary conditions. So $\vphi$ has two energy groups $\vphi=(\vphi^1,\vphi^2)$. Index 1 denotes the high energy group and 2 the thermal energy one. The flux is the solution to the following eigenvalue problem (see \cite{Hebert2009})
\begin{center}
Find $(\lambda,\phi) \in \bC \times \left(H^1(\cR)\times H^1(\cR)\right)$ such that for all $x\in \cR$,
\end{center}
\begin{equation}
\label{eq:diffusion}
\begin{cases}
&-\nabla\left(D^1(r)\nabla\vphi^1(r)\right)+\Sigma_{a}^1(r)\vphi^1(r) - \Sigma_{s,0}^{2\to1}(r)\vphi^2(r)= \lambda \left( \chi^1(r)(\nu\Sigma_{f})^1(r)\vphi^1+ \chi^1(r)(\nu\Sigma_{f})^2(r)\vphi^2(r) \right) \\
&-\nabla\left(D^2(r)\nabla\vphi^2(r)\right)+\Sigma_{a}^2(r)\vphi^2(r)-\Sigma_{s,0  }^{1\to2}(r)\vphi^1(r)
=\lambda\left( \chi^2(r)(\nu\Sigma_{f})^1(r)\vphi^1(r)+ \chi^2(r)(\nu\Sigma_{f})^2(r)\vphi^2(r) \right),
\end{cases}
\end{equation}
with
$$D^g(r)\nabla \phi^g(r)\cdot n + \frac{1}{2}\phi^g(r)=0 \quad  \text{ on }\partial \cR, \text{ for }g=1,\,2.
$$
The coefficients involved are the following:
\begin{itemize}
\item $D^g(r)$ is the diffusion coefficient of group $g$ with $g\in\{1,2\}$.
\item $\Sigma_{a}^g(r)$ is the macroscopic absorption cross section of group $g$.
\item $\Sigma_{s,0}^{g'\to g}(r)$ is the macroscopic scattering cross section of anisotropy order 0 from group $g'$ to $g$.
\item $\chi^g(r)$ is the fission spectrum of group $g$. 
\end{itemize}
We assume that they are either constant of piecewise constant in $\cR$ so we can view them as functions from $L^\infty(\cR)$.

The generated power is
\begin{equation}
\label{eq:power}
P\coloneqq (\kappa\Sigma_{f})^1 \phi^1 + (\kappa\Sigma_{f})^2 \phi^2,
\end{equation}
and since $\phi^1$ and $\phi^2\in H^1(\cR)$, we have $P\in V$.
 
We next make some comments on the coefficients and recall well-posedness results of the eigenvalue problem \eqref{eq:diffusion}. First of all, the first four coefficients ($D^g$, $\Sigma_{a}^g$, $\Sigma_{s,0}^{1\to2}$, $\Sigma_{s,0}^{2\to1}$  and $(\nu\Sigma_{f})^g$) might depend on the spatial variable. In the following, we assume that they are either constant or piecewise constant so that our set of parameters is
\begin{equation}
\label{eq:paramDiff}
\mu
=
\{
D^1, D^2, \Sigma_{a}^1, \Sigma_{a}^2, \Sigma_{s,0}^{1\to2}, (\nu\Sigma_f)^1, (\nu\Sigma_{f})^2, \chi^1, \chi^2
\}.
\end{equation}
By abuse of notation, in \eqref{eq:paramDiff} we have written $D^g$ to denote the set of values that this coefficient might take in space and similarly for the other parameters.

Under some mild conditions on the parameters $\mu$, the eigenvalue $\lambda_{\min}$ with the smallest modulus is simple, real and strictly positive (see \cite[Chapter XXI]{DLvol6}). The associated eigenfunction $\phi$ is also real and positive at almost every point $x\in\cR$ and it is what is classically called the flux. In neutronics, it is customary to work with the inverse of $\lambda_{\min}$, which is called the multiplication factor
\begin{equation}
\keff\coloneqq 1/ \lambda_{\min}.
\end{equation}
Therefore $\keff$ is not a parameter in our setting because, for each value of the parameters $\mu$, $\keff$ is determined by the solution to the eigenvalue problem. 

If the parameters of our diffusion model  range in, say, $$D^1 \in [D^1_{\min}, D_{\max}^1],\ D^2 \in [D^2_{\min}, D_{\max}^2],\dots,\chi^2\in[\chi^2_{\min}, \chi_{\max}^2],$$ then
\begin{equation}
\cD 
\coloneqq
[D^1_{\min}, D^1_{\max}]\times\dots\times[\chi^2_{\min}, \chi^2_{\max}],
\end{equation}
and the set of all possible states of the power is given by
\begin{equation}
\label{eq:manifold_power}
\cM_{\diffusion}
=
\{
P(\mu)\ :\ \mu \in \cD
\}\; \subset V,
\end{equation}
which is the manifold of solutions of our problem. 

\section{Numerical Examples}
\label{sec:numerics}

\subsection{Description of the test case and the numerical solver}

\medskip

\begin{figure}[htbp]
\begin{center}
\begin{tikzpicture}[scale=5./25.]
\fill[white] (0,0) rectangle (25.,25.);
\fill[red] (0,0) rectangle (15.,15.);
\fill[blue] (15.,0.) rectangle (20.,5.);
\draw (0,0) grid[xstep=5.,ystep=5.](25.,25.);
\draw[black](26.38,23.17) rectangle (29.512,24.2159);
\fill[white] (26.38,23.17) rectangle (29.512,24.2159);
\node at (33.953,23.851925) {Reflector };
\draw[black](26.38,21.17-0.1) rectangle (29.512,22.2159-0.1);
\fill[red] (26.38,21.17-0.1) rectangle (29.512,22.2159-0.1);
\node at (33.953,21.851925-0.1) {Core };
\draw[black](26.38,19.17-0.1) rectangle (29.512,20.2159-0.1);
\fill[blue] (26.38,19.17-0.1) rectangle (29.512,20.2159-0.1);
\node at (33.953,19.851925-0.1) {Void };
\draw[black](26.,18.17) rectangle (37.512,25.);
\node at(0.5*25.,-0.1*25.) { Reflexion};
\node at(-0.15*25.,0.5*25.) { Reflexion};
\node at(0.5*25.,+1.1*25.) { Vacuum};
\node at(+1.15*25.,0.5*25.) { Vacuum};
\end{tikzpicture}
\end{center}
\caption{Cross-sectional view of the core ($z = 0$ cm).}
\label{fig:takeda1}
\end{figure}
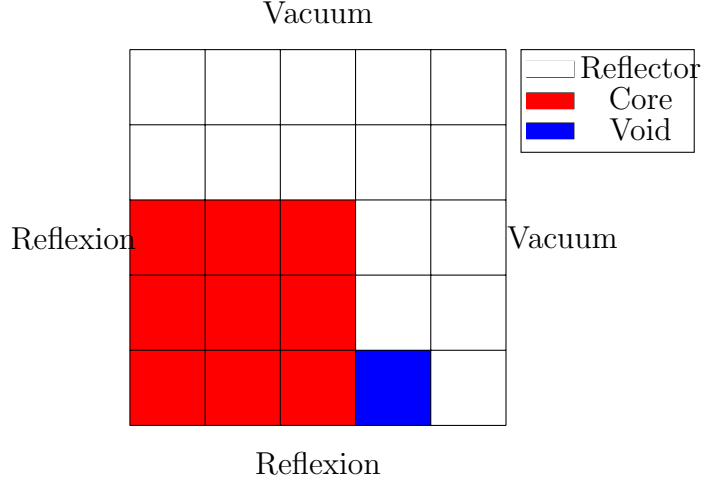

\textbf{The test-case:} We consider Model 1 Case 1 of the well-known Takeda neutronics benchmark~\cite{takeda1991benchmark} to build our test case. The
geometry of the core is three-dimensional and the domain is $\{(x,y,z)\in\mathbb{R}^3, 0\leq x\leq 25\,\text{cm};0\leq y\leq 25\,\text{cm};0\leq z\leq 25\,\text{cm}\} $. This test is defined with $G=2$ energy groups and isotropic scattering and we set $\kappa^g = 1$ MeV for $g=1,2$. The reactor core geometry is depicted in Figure \ref{fig:takeda1}. 
In the following, we implicitly refer to the cross-sections and the other coefficients of this test case.
Our goal is to reconstruct in real time the spatial power field of the reactor. We assume that the neutron transport equation perfectly describes reality, and the set of all possible states is given by the manifold
$$
\cM_{\transport}
=
\{
P^{\transport}(\mu)\ :\ \mu \in \cD
\}\; \subset V.
$$
The set of solutions of the neutron diffusion equation is
$$
\cM_{\diffusion} 
=
\{
P^{\diffusion}(\mu)\ :\ \mu \in \cD
\}\; \subset V.
$$
It is an imperfect description of the true states given by $\cM_{\transport}$.

The parameter set $\mu$ from equation \eqref{eq:paramDiff} is generated by the mapping
\begin{align*}
\mu : \left[0.8,1\right]^5 \subset \bR^5 &\to \bR^9 \\
\; \alpha &\mapsto \mu(\alpha)\left(
\frac{D^1}{\alpha_1}, \frac{D^2}{\alpha_2},\alpha_1\Sigma_{a}^1, \alpha_2\Sigma_{a}^2, \alpha_3\Sigma_{s,0}^{1\to2}, \alpha_4(\nu\Sigma_f)^1, \alpha_5(\nu\Sigma_{f})^2, \chi^1, \chi^2
\right).
\end{align*}
We can thus view the parameter set either as the $5$ dimensional tensorized subset $\left[0.8,1\right]^5$ where $\alpha$ ranges, or as a 5-dimensional surface manifold from $\bR^9$ where the $9$ coefficients $\mu$ of the neutronic model live.

We work with $m=54$ measurements observations that are placed uniformly in the reactor. They are defined as local averages over small subdomains $\cR_i\subset \cR$
\begin{equation}
    \label{eq:measures}
\omega_i(x) = \frac{1}{|\cR_i|} \mathds{1}_{\cR_i}(x),\quad \forall x\in \cR
,\,i=1,\dots,m.
\end{equation}
We compare two cases:
\begin{enumerate}
\item \textbf{Perfect physical model:} We apply PBDW using reduced models from the transport manifold which represents the true reality in our experiments.
\item \textbf{Imperfect physical model:} We assume that a perfect model is out of reach and we use the diffusion manifold. The reconstruction will thus be affected by a model bias.
\end{enumerate}

\medskip

\textbf{The solver:} To generate the snapshots and the reduced models, we have worked with MINARET~\cite{lautard_minaret}, a deterministic solver for reactor physics calculations developed in the framework of the APOLLO3$^{\text{\textregistered}}$\ code~\cite{schneider_apollo3}. MINARET can solve either the multigroup neutron transport or diffusion problem from Equations \eqref{eq:boltzmann_void} and \eqref{eq:diffusion}. The numerical scheme to compute the multiplication factor $\keff$ is based on the inverse power method (see, e.g., \cite{Hebert2009}). MINARET uses the $S_N$ discrete ordinate method to deal with the angular variable, and Discontinuous Galerkin Finite Elements to solve spatially the neutron transport equation~\cite{reed_triangular}. It applies the Symmetric Interior Penalty Galerkin method (SIPG)~\cite[Chapter 4]{dipietro2012mathematical} for the discretization of the neutron diffusion equation~\eqref{eq:diffusion}. In all cases, the solver uses cylindrical meshes devised by extrusion of a 2D triangular mesh. For our simulations, we work with a level-symmetric formula of order $N=8$ for the S$_N$ quadrature, and the spatial approximation uses discontinuous $\bP_1$ finite elements of a uniform mesh. 
The physical output power map is post-processed on an approximation space of dimension $N_h=540$ ($N_h$ degrees of freedom). 

\subsection{Case 1: Reconstruction with a perfect physical model}
\label{transport-results}

Here we assume that we have access to a perfect description of the physics, and we work with the neutron transport manifold $\cM_{\transport}$. 

In order to create a reduced space $V_n$ of small dimension $n \ll N_h$, we apply a Proper Orthogonal Decomposition (POD) based on the training set
\begin{align*}
\mathcal{P}_{\text{training}} =\{P^{\transport}(\mu(\alpha)), \alpha \in \{0.8,0.9,1\}^5\}\subset \cM_{\transport}
\end{align*}
of 
power maps obtained from solutions of the transport neutron equations, also called snapshots. 

\begin{figure}[h] 
\centering
\includegraphics[scale=0.6]{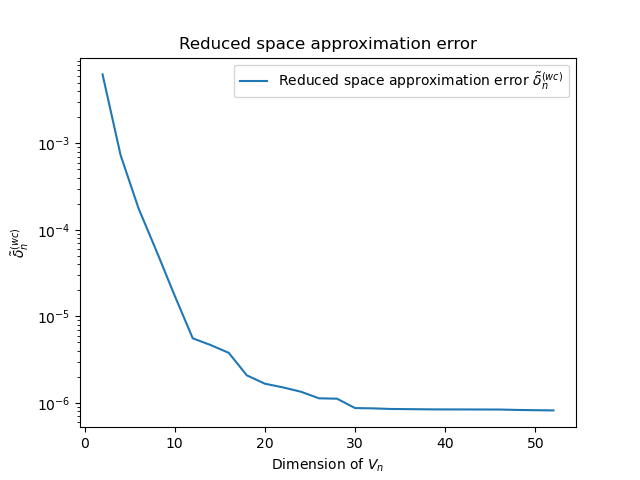}
\caption{Relative approximation error $\widetilde{\delta}_n^{(\wc)}$ of the transport manifold $\cM_{\text{tr}}$ with respect to the dimension $n$ of the reduced space. Here the reduced space is a POD computed using the same manifold $\cM_{\text{tr}}$.}
\label{fig:eps_N_transp}
\end{figure}

We measure the relative approximation error $\widetilde{\delta}_n^{(\wc)}$ as defined in Equation \eqref{eq:error-manifold}. For this, we define a collection of power maps of reference
\begin{equation}
\mathcal{P}_{\text{test}} = \{
P^{\transport}(\mu(\alpha))\ ,\ \alpha \in \{0.85,0.95\}^5
\}.
\label{eq:Ptest}
\end{equation}
Figure \ref{fig:eps_N_transp} shows that the training space is well approximated with a few POD modes. For $n \geq 30$, the relative error between one power map and its projection onto $V_n$ is smaller than $10^{-6}$. 

\begin{figure}[h] 
    \centering
    \includegraphics[scale=0.6]{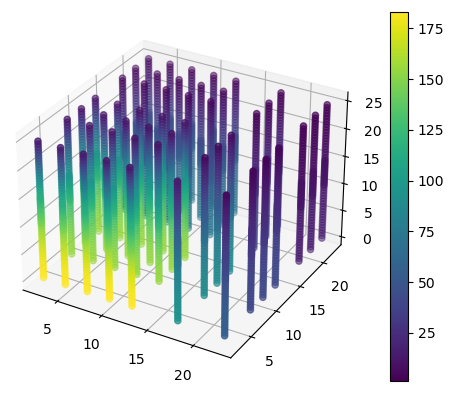}
    \includegraphics[scale=0.6]{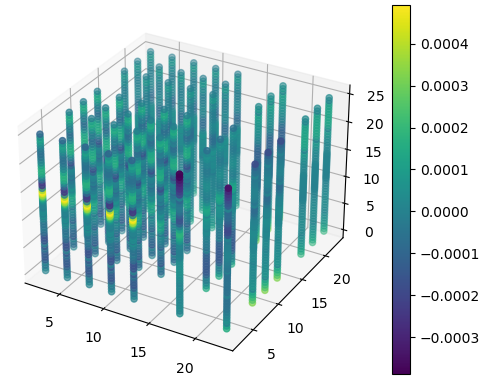}
    \includegraphics[scale=0.6]{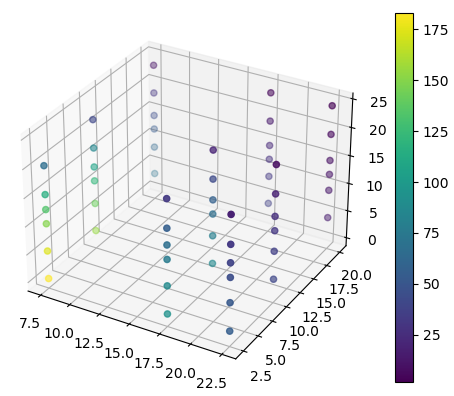}
    \caption{3D representation of the power map $P^{\transport}(\mu(\alpha))$ with $\alpha=\{0.85\}^5$ (upper left), the algebraic reconstruction error by PBDW (upper right) and the $m=54$ measurements}
    \label{fig:plots_3D}
\end{figure}

We next study the ability to reconstruct the power field with measurement observations, and the PBDW method, as Figure \ref{fig:plots_3D} shows in the 3D space. For this, we compute for $1\leq n \leq m$:
\begin{itemize}
\item The relative reconstruction error given by $\widetilde{e}_{m,n}^\wcpbdw = \max_{u \in \widetilde \cM_\transport}\dfrac{\Vert u - A_{m,n}^{\pbdw}(\omega)\Vert}{\Vert u \Vert}$,
\item The upper bound of the reconstruction error given by $\beta^{-1}(V_n, W_m)\widetilde{\delta}_n^{(\wc)}$, as given in Equation \eqref{eq:err-ms-pbdw}.
\end{itemize}

\begin{figure}[h] 
    \centering
    \includegraphics[scale=0.6]{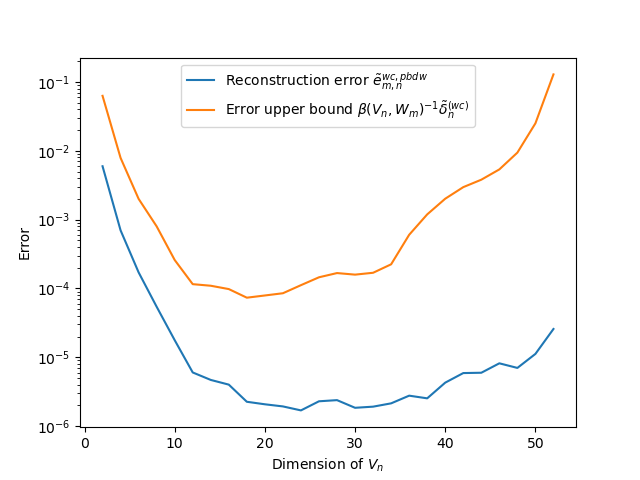}
    \caption{Relative reconstruction error $\widetilde{e}_{m,n}^\wcpbdw$ (in blue) and error estimate $\beta(V_n,W_m)^{-1}\widetilde{\delta}_n^{(\wc)}$ (in yellow) with respect to the dimension $n$ of the reduced space}
    \label{fig:rec_error_transp}
\end{figure}

Figure \ref{fig:rec_error_transp} shows 
that the upper bound is about two orders of magnitude above the actual reconstruction error. This gap is expected to decrease if we use more  functions in the test set. The second observation is that 
the reconstruction accuracy reaches a minimum for a dimension $n^*\approx25$. If we work with the optimal dimension $n^*$, an important result is that we can recover the power field from measurement observations at almost the same accuracy ($\approx 10^{-6}$, see Figure \ref{fig:eps_N_transp}) as the one given by the orthogonal projection onto $V_n$ (to see this, compare the errors at $n^*$ in Figures \ref{fig:eps_N_transp} and \ref{fig:rec_error_transp}).

\begin{figure}[h]
\centering
\includegraphics[scale=0.6]{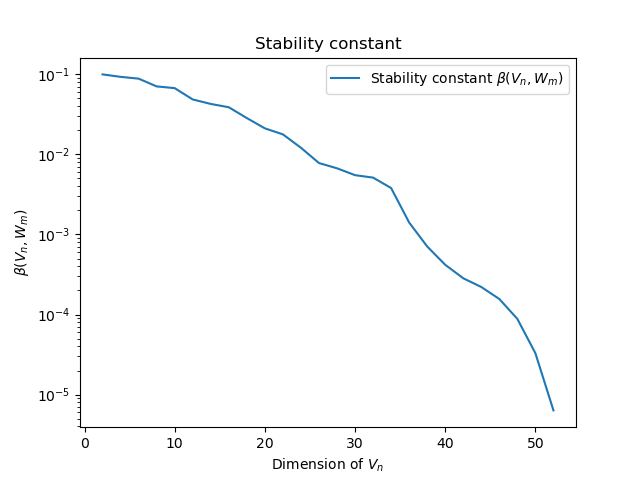}
\caption{Stability constant $\beta(V_n, \Wm)$ with respect to the dimension $n$ of the reduced space}
\label{fig:beta_transp}
\end{figure}

The behavior of the reconstruction error with the dimension $n$ is connected to a loss of stability illustrated in  Figure \ref{fig:beta_transp}.
It warns about a compromise to find between the approximation error of the manifold and the stability 
 in order to optimize the accuracy of the power map reconstruction. One strategy to mitigate stability problems is to find locations for the sensor measurements that span spaces $W_m$ maximizing the value of $\beta(V_n, W_m)$ (see, e.g., \cite{BCMN2018}).


\subsection{Case 2: Reconstruction of the power map from diffusion snapshots}
\label{diffusion-results}

We now consider the diffusion neutron equations as the best available model while the true states are given by the neutron transport model. They are therefore members of $\cM_{\text{tr}}$.

Similarly as done in Section \ref{transport-results}, we apply a POD over a collection
\begin{align*}
\mathcal{P}_{\text{training}} =\{P^{\diffusion}(\mu(\alpha)) \cond \alpha \in \{0.8,0.9,1\}^5\}\subset \cM_{\diffusion},
\end{align*} 
  to create a reduced space $V_n$ of dimension $n \ll N_h$. 
The main difference lies in the fact that the snapshots are obtained from the neutron diffusion equations, as we consider that the transport model cannot be computed in this section.

\begin{figure}[h]
\centering
\includegraphics[scale=0.6]{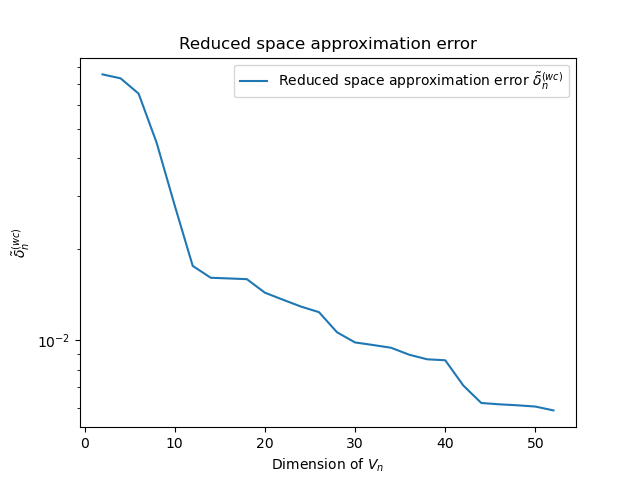}
\caption{Relative approximation error $\widetilde{\delta}_n^{(\wc)}$ of the transport manifold $\cM_{\text{tr}}$ with respect to the dimension $n$ of the reduced space. Here the reduced space is a POD computed using the diffusion manifold $\cM_{\text{diff}}$.}
\label{fig:eps_N_diff}
\end{figure}

Figure \ref{fig:eps_N_diff} shows that the approximation error of the transport manifold $\cM_{\text{tr}}$ is less accurate than in the previous case
due to the bias between the two models. Typically, for $n=50$, we approximate the manifold at the accuracy of $6 \times 10^{-3}$, whereas the approximation with the transport model was about $10^3$ times better (compare Figure \ref{fig:eps_N_diff} and Figure \ref{fig:rec_error_transp}). Therefore, the reconstruction error will have a similar order of magnitude to those observed for the approximation error.  

Similarly, we compute for $1\leq n \leq m$:
\begin{itemize}
\item The relative reconstruction error given by $\widetilde{e}_{m,n}^\wcpbdw = \max_{u \in \widetilde \cM_{\transport}}\dfrac{\Vert u - A_{m,n}^{\pbdw}(\omega)\Vert}{\Vert u \Vert}$,
\item The upper bound of the reconstruction error given by $\beta^{-1}(V_n, W_m)\widetilde{\delta}_n^{(\wc)}$, as given in Equation \eqref{eq:err-ms-pbdw}.\\
\end{itemize}

\begin{figure}[h]
\centering
\includegraphics[scale=0.6]{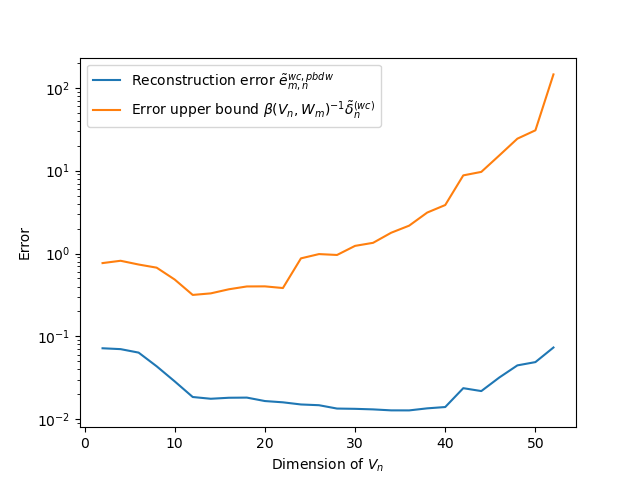}
\caption{Relative reconstruction error $\widetilde{e}_{m,n}^\wcpbdw$ (in blue) and error estimate $\beta(V_n,W_m)^{-1}\widetilde{\delta}_n^{(\wc)}$ (in yellow) with respect to the dimension $n$ of the reduced space}
\label{fig:rec_error_diff}
\end{figure}

As done before, the PBDW reconstruction procedure is then performed by extracting measurements over the collection power maps of reference defined in~\eqref{eq:Ptest}. Figure \ref{fig:rec_error_diff} illustrates that
 the minimum for the reconstruction error reaches about $1.5 \times 10^{-2}$ for $n^* \approx 35$. The gap between the reconstruction error and its estimate here is bigger as the stability plays a secondary role. Hence, the reconstruction error is only led by the model bias.

\begin{figure}[h]
\centering
\includegraphics[scale=0.6]{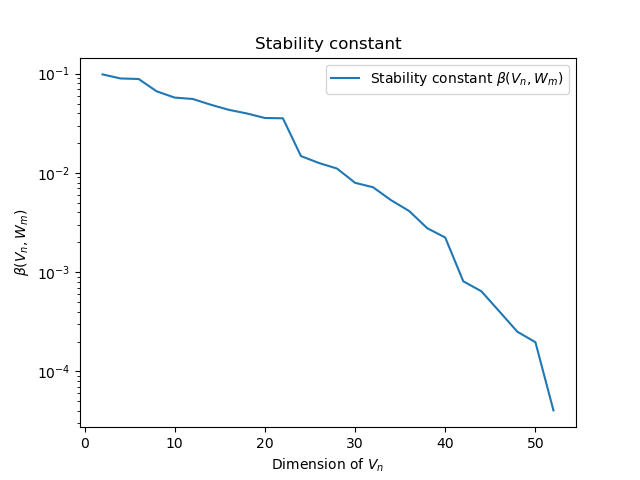}
\caption{Stability constant $\beta(V_n, \Wm)$ with respect to the dimension $n$ of the reduced space}
\label{fig:beta_diff}
\end{figure}

Figures \ref{fig:beta_diff} shows that the stability constant
 presents the same behavior as in the case of $V_n$ built with transport snapshots. 

\newpage

\section*{Conclusion of the numerical results and outlook}

From the numerical experiments, it follows that the PBDW algorithm can reconstruct data very efficiently when the physical model is perfect. An interesting fact is that there are optimal values $n^*$ for the dimension in the reduced models $V_n$ used in the PBDW algorithm which make the reconstruction with measurement observations be comparable to the approximation accuracy by projection on $V_n$ (see, e.g., Figures \ref{fig:eps_N_transp} and \ref{fig:rec_error_transp}).

In presence of model error, the conclusions are analogous. However, the approximation accuracy by direct projection is degraded by the presence of the model error as the comparison between Figure \ref{fig:eps_N_transp} and Figure \ref{fig:eps_N_diff} shows. \textcolor{black}{This degradation may be reduced if some snapshots are computed with the transport model. The selection of these snapshots may be based on a posteriori estimators devised specifically for the model error. }

In principle, the PBDW is expected to be able to correct to some extent the model error due to fact that reconstructions lie in $V_n\oplus (W\cap V_n^\perp)$ and not only in $V_n$. There, if the model is biased and yields a reduced model $V_n$ which is not perfectly appropriate, the component $(W\cap V_n^\perp)$ is expected to help to correct this inaccuracy. However, our results tend to indicate that this correction component has a very limited effect in our case. This may be due to the poor approximation properties of the observation space $W$, which is, in our case, spanned by functions that are very localized in space (see equation \eqref{eq:measures} for the definition of the $\omega_i$). This behavior could be improved by working with parametrized families of spaces such as Reproducing Kernel Hilbert spaces (see \cite{maday2019adaptive}).
 In that case, we could try to find an appropriate space for which the $\omega_i$ would better enhance the final reconstruction quality. Another option would be to consider purely data-driven corrections on top of the PBDW reconstruction, making use of supervised learning techniques and feed forward neural networks. These ideas will be the starting point of future works in mitigating the effect model error in state estimation.

%

\section*{Acknowledgements}
The authors gratefully acknowledge O. Lafitte for fruitful discussions.
This project was partly funded by CEA. O.M.~was funded by the Emergence project grant of the Paris City Council ``Models and Measures''.

\bibliographystyle{unsrt}
\bibliography{literature}
\end{document}